# LINEAR HYPERGRAPH EDGE COLORING

*Vance Faber*
*Revision: March 3, 2016*

**Abstract**. Motivated by Erdos-Faber-Lovász (EFL) conjecture for hypergraphs, we consider the edge coloring of linear hypergraphs.

**The back story**. In 1972, a three week conference on hypergraphs was held at The Ohio State University. It was during this conference that the seeds of EFL were planted. In particular, in Problem 18 (see below and [7]) we asked for bounds on edge coloring. Later, we added the condition that the hypergraph be linear which modified the conjectured bounds in Problem 18 and when we couldn't prove or disprove that we added the additional constraints that created EFL. In this paper, I go back to the precursor of EFL and show that some of the facts we know about EFL apply equally well to the precursor.

**Notation**. Let $H = (V, E)$ be a *hypergraph* (see, for example [6]): a set of subsets $E$ of the set $V$. We call the elements of $V$ the *vertices* and the elements of $E$ the *edges*. We often write $n = |V|$ and $m = |E|$. The *degree* of a vertex $x$ is the number of edges $d(x)$ which include it. We let the minimum degree be $\delta$ and the maximum degree be $\Delta$. If all vertices have the same degree, we say the hypergraph is *regular*. The *rank* of an edge $e$ is the cardinality $r(e)$ of $e$. We let the minimum rank be $\rho$ and the maximum rank be $P$. If all edges have the same rank, we say the hypergraph is *uniform*. If $r(e) = 2$ for every edge then $H$ is a *graph*. If the intersection of any two edges has at most one vertex, we call the hypergraph *linear*.

**Incidence matrix formulation**. An equivalent formulation for a hypergraph is to consider $H$ to be the incidence matrix of the hypergaph. In this case, $H$ is an $n \times m$ matrix: a row of $H$ is the transpose of the characteristic vector of a vertex and a column of $H$ is the characteristic vector of a edge. We use these two formulations interchangeably. It is often easier to understand a fact in one formulation or the other. For example, a fundamental theorem for a hypergraph is that the sum of the ranks is equal to the sum of the degrees. This is trivial to see in the matrix formulation because both sides of the equality are clearly equal to the number of non-zero entries in the matrix $H$. In this formulation, an edge $e$ is a column vector and a vertex $x$ is a row vector. Two vertices $x$ and $y$ are independent if and only if they are orthogonal, that is, the inner product $xy^* = (x^*, y^*) = 0$. Two edges are independent if and only if they are orthogonal, that is $e^* f = (e, f) = 0$.

**Duality**. If $H$ is a hypergraph, so is the transpose $H^*$ called the *dual* hypergraph. The theorems we need are often stated in a form that applies most naturally to $H^*$ and we have to translate them to $H$ to apply them. Clearly the edges of $H^*$ are the vertices of $H$ and vice versa, ranks swap with degrees, etc. It is easy to see that $H$ is linear if and only if $H^*$ is linear. This is because $H$ is linear if and only if every $2 \times 2$ minor has a zero entry.

**Clique and line graphs**. The clique graph $C(H)$ has the same vertex set as $H$ and an edge for every pair of vertices in some edge. Each edge in $H$ then appears as a clique in $C(H)$. The line graph of $H$ is the clique graph of $H^*$; it has a vertex for each edge of $H$ and an edge between two edges of $H$ if they intersect. Note that in a linear hypergraph, the cliques in the clique graph which come from the edges in the hypergraph are edge disjoint; the clique graph is a set of edge disjoint complete subgraphs of the complete graph on $n$ vertices.

**Coloring**. A coloring of the edges of a hypergraph is a function $\gamma$ from the edges of a hypergraph into a set $\Gamma$, called colors, such that $\gamma(e) = \gamma(f)$ only if $e$ and $f$ are orthogonal. We let $q(H)$, called the chromatic index of $H$, be the cardinality of the smallest $\Gamma$ for which there is a coloring. The chromatic index of $H$ is the chromatic number of the line graph. If $H$ is colored with $q$ colors, the columns of the matrix $H$ are the disjoint union of $q$ sets of columns, each one of which consists of pairwise orthogonal columns.

**Conjectures.** There is a long standing conjecture known now as the Erdos-Faber-Lovász conjecture which in its edge formulation says

*Conjecture C1 (EFL). Let $H$ be a linear hypergraph with $n$ vertices and no rank 1 edges. Then $q(H) \leq n$.*

In fact, there are stronger conjectures from which EFL stemmed that appear to have been forgotten.

*Definition.* We define the *clique degree* of a vertex $x$ by

$$D(x) = \sum_{e \supset x} (r(e) - 1).$$

This is the degree of $x$ in the clique graph. The dual concept is the *clique rank* of an edge

$$R(e) = \sum_{x \in e} (d(x) - 1).$$

*Conjecture C2. Let $H$ be a linear hypergraph with maximum rank $\mathrm{P}$, maximum degree $\Delta$ and no rank 1 edges such that for every vertex $D(x) \leq k$. Then $q(H) \leq k + 1$.*

*Conjecture C3. Let $H$ be a linear hypergraph with maximum rank $\mathrm{P}$, maximum degree $\Delta$ and no rank 1 edges then $q(H) \leq \Delta \mathrm{P} - \max(\Delta, \mathrm{P}) + 1$.*

*Discussion about designs.* Most examples that we use come from designs. A $(v, r, 1)$ *balanced incomplete block design* ($(v, r, 1)$ *BIBD*) is a linear hypergraph with $v$ vertices, rank $r$ and every pair of vertices is in exactly one edge. A *projective $(n, r, \Delta, 1)$ design* (see [5]) is a uniform regular hypergraph with $n$ vertices, rank $r$ and degree $\Delta$ and every two edges have a non-empty intersection. A projective design is not required to be linear. Note that the dual of a $(v, k, 1)$ BIBD is projective design with parameters $\Delta = k$, $r = (v-1)/(k-1)$, and $n = (v/k)(v-1)/(k-1)$. We know that if $H$ is a projective design, then $q(H)$ is the number of edges but we often do not know what $q(H^*)$ is. A *projective plane of order $r$* is a projective $(r^2 + r + 1, r + 1, 1)$ BIBD; its dual is also a projective design.

*Conjecture C4. Let $H$ be a linear hypergraph with maximum rank $\mathrm{P} \geq 3$, maximum degree $\Delta$ and no rank 1 edges then $q(H) \leq \Delta \mathrm{P} - \max(\Delta, \mathrm{P})$ unless $\Delta \leq \mathrm{P}$ and $H$ is a projective design in which case $q(H) \leq \Delta \mathrm{P} - \mathrm{P} + 1$.*

*Conjecture C5. Let $H$ be a linear hypergraph with maximum rank $\mathrm{P} \geq 3$, maximum degree $\Delta$ and no rank 1 edges such that for every vertex $D(x) \leq k$. Then $q(H) \leq k$ unless $\Delta = \mathrm{P} = \mathrm{T}$, $k = \mathrm{T}^2 - \mathrm{T}$ and $H$ is a $(\mathrm{T}^2 - \mathrm{T} + 1, \mathrm{T}, 1)$ BIBD in which case $q(H) = k + 1$.*

**Background Facts (edge formulations).**

*Fact 1 (Brooks' Theorem - edge formulation). For any hypergraph if $R(e) \leq R$ for every edge $e$ then $q(H) \leq R$ unless the clique graph of the dual is either an odd cycle or a complete graph in which case the bound is $R + 1$.*

Note that if the clique graph of the dual is an odd cycle, then the dual is an odd cycle also so $R = 2$ while $q = 3$. If the clique graph is a complete graph, then every pair of vertices in $H^*$ is in some edge so every pair of edges in $H$ meets and $q(H) = R + 1$.

*Fact 2 (Vizing's Theorem). For any graph $G$, if $d(x) \leq k$ then $q(H) \leq k+1$.*

Note that Conjecture C2 is a natural extension of Vizing's theorem. If you think of $H$ as an edge disjoint union of cliques, then the hypothesis of Conjecture C2 is that the sum total of all the edges of all the cliques that meet at a fixed vertex $x$ is at most $k$. Vizing's theorem colors these edges so that they all have different colors. Conjecture C2 demands that all the edges that belong to the same clique have the same color.

*Definition. Three edges $(e_1, e_2, e_3)$ in $H$ are called a triangle with side $e_1$ if either all the edges have a non-empty intersection or they have pair-wise non-empty intersections. We denote by $T(e)$ the number of triangles in $H$ with side $e$.*

*Fact 3 (AKS Theorem – edge formulation). Let $H$ be a hypergraph with $R(e) \leq R$. Suppose for every edge, $T(e) \leq R^2/f$. Then there exists a universal constant $c$ such that $q(H) \leq cR/\log f$.*

**Some theorems.** In this section we use the facts to prove theorems.

*Theorem 4.* $C5 \Rightarrow C4 \Rightarrow C3$, $C5 \Rightarrow C2 \Rightarrow C3$ and $C2 \Rightarrow C1$.

*Proof.* Assuming C5, let $H$ be a linear hypergraph with $P \geq 3$, maximum degree $\Delta$ and no rank 1 edges. We have two cases depending upon whether $P > \Delta$ or not. If $P > \Delta$, then we employ Brooks' Theorem so that we have $q(H) \leq \Delta P - P$ unless every two edges have a non-empty intersection. By counting the edges that meet a given edge $e$, we find that

$$q(H) \leq \Delta P - P + 1$$

with equality only if $e$ has rank $P$ and every vertex on $e$ has degree $\Delta$. But then $H$ is a projective design. If $P \leq \Delta$ then by C5 $q(H) \leq \max D(x) \leq \Delta P - \Delta$ unless $P = \Delta$ and $H$ is a BIBD. By counting the number of intersections, it is easy to see that $H$ is projective. This proves C4. To get C3 from C4, if $P = 2$ this is Vizing's Theorem. Otherwise, it follows directly from C4.

   Proving C2 from C5 is obvious. To prove C3 from C2, we copy the proof of obtaining C4 from C5 but there is no need to deal with the projective complication.

   Assuming C2, let $H$ be a linear hypergraph with $n$ vertices and no rank 1 edges and let $x$ be a vertex. Since $H$ is linear, all the vertices besides $x$ in all the edges that meet $x$ are distinct so there can be at most $n-1$ of them. Thus $D(x) \leq n-1$ and so by conjecture C2, $q(H) \leq n$. This proves C1.

*Theorem 5. Suppose that $H$ is a linear hypergraph of rank $\mathrm{P}$ with no rank 1 edges. There is a universal constant $C$ such that if $\mathrm{P} \geq C \geq 3$ and $\Delta \geq C(\mathrm{P}-1)$ then*

$$q(H) \leq \Delta(\mathrm{P}-1) - 1.$$

*Proof.* If $H$ is not uniform, add vertices of degree one as needed to make it uniform. Thus we can assume that $H$ is uniform with the rank of every edge equal to $r = \mathrm{P} \geq 3$. Let $k = \Delta(r-1)$ and assume that $H$ is not $k-1$ colorable. We work with the line graph $L = L(H)$ of the hypergraph $H$. The line graph uses the $m$ edges of $H$ as vertices and two edges in $H$ form an edge in $L$ if they meet in $H$. Any coloring of the vertices of $L$ is a coloring of the edges of $H$ and vice versa. We are given for every $x$

$$k \geq (r-1)d(x)$$

so

$$d(x) \leq \frac{k}{r-1}. \qquad (1)$$

Thus

$$d(x) - 1 \leq \frac{k}{r-1} - 1 = \frac{k+1-r}{r-1}.$$

Thus degree $R(e)$ of a vertex $e$ in the line graph $L$ can be no more than

$$R(e) = \sum_{x \in e}(d(x) - 1) \leq r\frac{k+1-r}{r-1}.$$

In particular, if we let

$$R = r\frac{k+1-r}{r-1} = (k+1)\frac{r}{r-1} - \frac{r^2}{r-1} \qquad (2)$$

the line graph has maximum degree at most $R$. Furthermore by (1) and (2),

$$d(x) \leq \frac{R}{r} + 1. \qquad (3)$$

By Brooks' Theorem, $L$ can be colored in $R$ colors unless it is a complete graph or an odd cycle. The odd cycle would give $q = 3$ so that is ruled out. If $L$ is complete, then $q = m$, the number of edges in $H$. Then for $H$ to be a counterexample, we would have $m > \Delta$. Consider a vertex $x$ of degree $\Delta$.

Edges not containing $x$ must meet every one of the $\Delta$ edges containing $x$ so they must have rank at least $\Delta$. But this contradicts $\Delta \geq C(r-1) \geq 3(r-1)$ since $r \geq 3$. Thus we only have to deal with cases where $R > k$. So we assume that
$$k \leq R. \tag{4}$$

We proceed by bounding the number $T(e)$ of triangles in $H$ with side $e$. This number is made up of two types of triangles which have edges $e$ and two other edges. First, there is $T_1(e)$, the triangles where all 3 edges have a common intersection. Second, there is $T_2(e)$, the triangles which meet pairwise at distinct vertices. The number of pairs of edges that meet a vertex $x$ in $e$ is given using (3) by

$$T_1(e) = \sum_{x \in e} \binom{d(x)-1}{2} \leq \frac{r}{2} \frac{R}{r}\left(\frac{R}{r}-1\right) = \frac{R}{2}\left(\frac{R}{r}-1\right).$$

To bound $T_2(e)$, we first select one of the $R(e)$ edges $e'$ that meets $e$ at some point $x$ and then select a pair of vertices $y \in e$ and $z \in e'$ to form the triangle with $x$. Since $H$ is linear these vertices can be in only one edge. This will count every triangle that exists twice so

$$T_2(e) \leq R \frac{(r-1)^2}{2}.$$

Thus

$$T(e) = T_1(e) + T_2(e) \leq \frac{R^2}{2}\left(\frac{1}{r} + \frac{(r-1)^2}{R}\right).$$

Let

$$\frac{1}{f} = \frac{1}{2r} + \frac{(r-1)^2}{2R}.$$

Then $T(e) \leq \dfrac{R^2}{f}$ so $q(H) \leq c \dfrac{R}{\log f}$ by AKS. Let $C = e^{2c}$. Since $r \geq C$, $2 - C/r \geq 1$. Since $\Delta \geq C(r-1)$, utilizing inequality (4) gives

$$R \geq k = \Delta(r-1) \geq C(r-1)^2 \geq \frac{C(r-1)^2}{2 - C/r}$$

or

$$\frac{(r-1)^2}{2R} \leq \frac{2 - C/r}{2C} = \frac{1}{C} - \frac{1}{2r}.$$

Thus

$$\frac{1}{f} = \frac{1}{2r} + \frac{(r-1)^2}{2R} \leq \frac{1}{C}.$$

Thus we have shown that under the given hypothesis that

$$q(H) \leq c \frac{R}{\log f} \leq \frac{R}{2}.$$

But from (2)

$$\frac{R}{2} = \frac{kr}{2(r-1)} + \frac{r}{2(r-1)} - \frac{r^2}{2(r-1)} \leq \frac{3}{4}k - \frac{r}{2} \leq \Delta(r-1) - 1.$$

*Corollary 6. Let $C$ be the universal constant from Theorem 5. Let $H$ be a hypergraph with $n$ vertices, $P \geq C \geq 3$, $\Delta \geq C(P-1)$ and no rank 1 edges. Then if $H$ is uniform or $n \geq \Delta(P-1) - 1$, $q(H) \leq n$.*

*Proof.* For every $x$, $D(x) \leq n-1$. If $H$ is uniform,

$$d(x)(P-1) = D(x) \leq n-1$$

and so

$$n > \Delta(P-1).$$

*Theorem 7 (see [3]). If $n > (\Delta-1)^2$ or $n < \rho^2$ then $q(H) \leq n$.*

*Proof.* Suppose we start with a minimal (with respect to the number of edges) counterexample $H$. Let the maximum degree of $H$ be $\Delta \geq 2$. If an edge $e$ has rank less than $n/(\Delta-1)$, then by induction we color $H \setminus \{e\}$ using $n$ colors and then since the number of colors sharing a vertex with $e$ is less than $n$, we can extend the coloring to all of $H$. Thus we can assume $\rho \geq n/(\Delta-1)$. Now let $x$ be a vertex of degree $\Delta$. Count the vertices adjacent to $x$. The number is at least

$$\Delta\left(\frac{n}{\Delta-1} - 1\right) = n + \frac{n}{\Delta-1} - \Delta$$

but can be no more than $n-1$ so

$$n + \frac{n}{\Delta-1} - \Delta \leq n-1.$$

Thus $n \leq (\Delta-1)^2$. Similarly, by counting the number of vertices adjacent to any vertex $x$ of degree $d(x)$,

so
$$d(x)(\rho-1) \leq n-1$$

$$\Delta \leq \frac{n-1}{\rho-1}.$$

Thus
$$n \leq \left(\frac{n-\rho}{\rho-1}\right)^2.$$

Solve this quadratic to get

$$n \geq \rho^2.$$

**List edge coloring.** We recall the definition of a list edge coloring of a graph. Suppose each edge of the graph $G$ has a set of colors (a list) associated with it. Then if the edges $G$ can be colored with colors chosen from the associated lists, then we say $G$ is *list colored*. If $G$ can be list colored using any arbitrary set of lists as long as they each have $k$ colors, then we say that $G$ is $k$ *list colorable*. We let $q_{list}$ be the smallest $k$ for which $G$ is $k$ list colorable.

*Conjecture C6 (Vizing's list coloring conjecture).* For every graph $G$, $q_{list}(G) \leq \Delta + 1$.

*Theorem 8.* Suppose C6 holds. Let $H$ be a linear hypergraph with $n$ vertices. Let $H_3$ be the hypergraph with edges of rank less than 3 removed. Let $H_2$ be the hypergraph with only the edges of rank 2. If the maximum degree of $H_3$ is $\Delta$ and the maximum degree of $H_2$, $\Delta(H_2)$, is at most $n - 2\Delta - 1$ then any $n$-coloring of the edges $H_3$ can be extended to an $n$-coloring of $H$.

*Proof.* The number of colors that are not used by the edges at $x$ in the coloring of $H_3$ is at least $n - \Delta$. Given an edge $e = (x, y)$ in $H_2$, the number of colors available to color $e$ is at least

$$(n-\Delta) + (n-\Delta) - n = n - 2\Delta \leq \Delta(H_2) + 1.$$

If Vizing's list edge-coloring conjecture holds, then each edge $e$ of $H_2$ can be colored using colors that are not used by other edges at the endpoints of $e$.

*Corollary 9.* Suppose C6 holds. If for every vertex $D(x, H_3) \geq 2\Delta$, then any $n$ edge coloring of $H_3$ can be extended to an $n$ edge coloring of $H$.

*Proof.* The degree of a vertex $x$ in the graph $H_2$ satisfies
$$d(x, H_2) = n - 1 - \sum_{e \supset x}(r(e) - 1) = n - 1 - D(x, H_3) \leq n - 1 - 2\Delta.$$
Thus $\Delta(H_2) \leq n - 1 - 2\Delta$.

*Corollary 10.* Suppose C6 holds. If $H_3$ is regular with degree $d$, then any $n$ edge coloring of $H_3$ can be extended to an $n$ edge coloring of $H$.

*Proof.* We can calculate $D(x, H_3) = \sum_{e \supset x}(r(e) - 1) \geq 2d$.

*Definition.* Let $d_k(x)$ be the number of edges with rank $k$ at the vertex $x$. Let $d(x)$ be the degree of the vertex $x$ in $H_3$. We call $\Delta - d(x)$ the deficit at $x$ and $\sum_{k \geq 4}(k - 3)d_k(x)$ the excess at $x$.

*Corollary 11.* Suppose C6 holds. If for every $x$ the excess at $x$ is at least as great as twice the deficit at $x$, then any $n$ coloring of $H_3$ can be extended to an $n$ coloring of $H$.

*Proof.* We have
$$d(x, H_2) = n - 1 - \sum_{e \supset x}(r(e) - 1) = n - 1 - 2d(x) - \sum_{k \geq 4}(k - 3)d_k(x) \leq n - 1 - 2\Delta.$$

**Critical hypergraphs.** Let $H$ be a linear hypergraph with maximum rank $P \geq 3$, maximum degree $\Delta$, no rank 1 edges and $H$ is not a $(P^2 - P, P, 1)$ BIBD. We let
$$D = D(H) = \max_x D(x) = \max_x \sum_{e \supset x}(r(e) - 1).$$
We want to investigate what properties $H$ must have to be a minimal counterexample to C2.

*Theorem 12.* If $H$ is a minimal counterexample to C2 then
(i) $q(H) > D(H)$;
(ii) for every $e$, $D(H \setminus e) = D(H)$;
(iii) for every $e$, $q(H \setminus e) = D(H)$;
(iv) for every $e$, $R(e) \geq D$.

*Proof.* (ii) If $D(H \setminus e) < D(H)$ then the minimality of $H$ implies that $q(H \setminus e) \leq D(H) - 1$. But then $q(H) \leq D$.

(iii) Consequently, if $q(H \setminus e) < D(H) = D(H \setminus e)$ then $q(H) \leq D$.

(iv) Suppose that $q(H) > D$ but for every edge $e$, $q(H \setminus e) \leq D$. If $e$ is an edge such that $R(e) < D$ then when the edges of $H \setminus e$ are colored in $D$ colors, only $D-1$ are used on the edges that meet $e$. Thus there is a color available for $e$ contradicting $q(H) > D$.

*Definition.* A *critical hypergraph* is one that satisfies the hypotheses of Theorem 12.

**Hypergraphs that are not necessarily linear**. We provide a partial solution to Problem 18 in Hypergraph Seminar.

*Problem 18 (Faber-Lovász).* Find the best possible function $f(\Delta, P)$ such that $q(H) \leq f(\Delta, P)$. It is easy to show

(i) $f(\Delta, P) \leq r(\Delta - 1) + 1$

(ii) $f(\Delta, P) = r(\Delta - 1) + 1$ if a $(k+1, \Delta, 1)$ BIBD exists with $k = r(\Delta - 1)$;

(iii) $f(\Delta, P) \geq \Delta(r^2 - 3r + 3)/(r - 1)$ when $r - 1$ divides $\Delta$.

*Discussion.* There may be an omission in part (iii). It certainly does not seem easy in general. The construction that I recall only works for $r$ a power of a prime plus 2. I summarize what I know about this problem in the following theorem.

*Theorem 13.* For every $r \geq 3$ and $\Delta \geq 2$, among the hypergraphs $H$ with maximum rank $r$ and maximum degree $\Delta$

*(i) there exists a uniform hypergraph with $q(H) = r(\Delta - 1) + 1$ if and only if a $(r(\Delta - 1) + 1, \Delta, 1)$ balanced incomplete block design exists;*

*(ii) there exists a uniform hypergraph with $q(H) = \Delta(r^2 - 3r + 3)/(r - 1)$ if a projective plane of order $r - 2$ exists (in particular, if $r$ is a power of a prime plus 2) and $r - 1$ divides $\Delta$.*
*(iii) there exists a uniform hypergraph with $q(H) = \Delta(r - 1)$ if a projective $(n, r, \Delta, 1)$ design exists with $n = r(r - 1)$.*

*Proof.* (i) Since we know the dual of the design has the right value for $q$, we only have to prove the necessity. We work with $C(H^*)$ which has maximum degree $k = r(\Delta - 1)$ and might be a multigraph if $H$ is not linear. We let $G$ be the graph formed by removing multiplicities from $C(H^*)$. Let $k = r(\Delta - 1)$. Suppose $q(H) = k + 1$. We know by Brooks' Theorem that $q(H) = k$ unless $G$ is a complete graph on $k+1$ vertices or an odd cycle. Since $r \geq 3$ a cycle is ruled out and so $G$ is a complete graph. But the vertices in $G$ are the vertices in $H^*$ so $H^*$ has $k+1$ vertices and every pair must be in at least one edge. Furthermore, the degree of each vertex in $G$ must be $k$. But the degree of the vertex in $G$ which corresponds to $e$ in $H$ is the sum over all vertices $x$ in $e$ of $d(x)-1$ minus any multiplicities. Suppose we have to remove $\gamma$ multiplicities. Then we have

$$r(\Delta - 1) - \gamma \geq \sum_{x \in e}(d(x) - 1) - \gamma = r(\Delta - 1)$$

so there can be no multiplicities and $G$ is regular with $d(x) = \Delta$. That means each edge in $C(H^*)$ is in exactly one clique, in other words, $H^*$ is linear.

(ii) Let $t = \Delta/(r-1)$. We duplicate each edge in the design $t$ times. Each edge in the design has cardinality $r-1$. We add to each edge $e$ a new vertex $x_e$. The new hypergraph now has edges of rank $r$, maximum degree $\Delta$ and

$$q = \frac{\Delta}{r-1}[(r-2)^2 + (r-2) + 1] = \frac{\Delta}{r-1}(r^2 - 3r + 3).$$

(iii) Suppose that $H$ is a projective $(n, r, \Delta, 1)$ design with $n = r(r-1)$. Since $H$ is projective, the number of edges must be $q$. Solving $qr = n\Delta = r(r-1)\Delta$ gives $q = \Delta(r-1)$.

*Remarks.* Note that this theorem leaves open the possibility that the extremal hypergraphs in Problem 18 are linear. Also, Lovász [4] showed that a projective $(n, r, \Delta, 1)$ design can exist only if $n \leq r^2 - r + 1$. The case $n = r^2 - r + 1$ is covered in part (i) of the theorem so $n = r^2 - r$ is the next smallest value.

**References**.

[1]  V. Faber, "The Erdős-Faber-Lovász conjecture – the uniform regular case", J. Combinatorics, 1(2010), 113-120.

[2]  N. Alon, M. Krivelevich and B. Sudakov, "Coloring graphs with sparse neighborhoods", *J. Comb. Theory*, B 77 (1999), 73-82.

[3]   A. Sánchez-Arroyo, "The Erdős-Faber-Lovász conjecture for dense hypergraphs", Discrete Mathematics 308(5-6) (2008), 991-992.

[4]  L. Lovász, "On minimax theorems of combinatorics (in Hungarian)," Mat. Lapok 26 (1975), 209-264.

[5]  P. Erdős, V. Faber and F. Jones, "Projective $(2n, n, \lambda, 1)$-designs", J. Stat. Plan. Inference 7 (1982), 181-191.

[6]  C. Berge, Hypergraphs: Combinatorics of Finite Sets, Vol. 45, North-Holland Mathematical Library, 1989.

[7]  Hypergraph Seminar (C. Berge and D. Ray-Chaudhuri, eds.), Vol 411, Lecture Notes in Mathematics, New York 1974.